\documentclass[12pt]{article}
\usepackage{tikz}
\usetikzlibrary{intersections,calc,arrows}
\usepackage{amssymb}
\date{}
\makeatletter
\newcommand{\figcaption}[1]{\def\@captype{figure}\caption{#1}}
\newcommand{\tblcaption}[1]{\def\@captype{table}\caption{#1}}

\@addtoreset{equation}{section}
\makeatother
\newcommand{\qed}{\hbox{\rule[-2pt]{3pt}{6pt}}}

\begin{document}
\title {\bf Asymptotic formulas for $L^2$ bifurcation curves of nonlocal logistic equation of population dynamics}

\author{{\bf Tetsutaro Shibata}
\\
{\small Hiroshima University, Higashi-Hiroshima, 739-8527, Japan}
}

\maketitle
\footnote[0]{E-mail: tshibata@hiroshima-u.ac.jp
}
\footnote[0]{Graduate School of Advanced Science and Engineering, 
Hiroshima University, Higashi-Hiroshima, 739-8527, Japan}

\footnote[0]{This work was supported by JSPS KAKENHI Grant Number JP25K07087.
}

\vspace{-0.5cm}

\begin{abstract}
The one-dimensional nonlocal Kirchhoff type bifurcation problems 
which are derived from logistic equation of population dynamics are 
studied.  
We obtain the precise asymptotic shapes of $L^2$ bifurcation curves 
$\lambda = \lambda(\alpha)$ as $\alpha \to \infty$, 
where $\alpha:= \Vert u_\lambda \Vert_2$.

\end{abstract}

\noindent
{{\bf Keywords:} Nonlocal Kirchhoff type elliptic equations, logistic equation of population dynamics, $L^2$-bifurcation curve} 

\vspace{0.5cm}

\noindent
{{\bf 2020 Mathematics Subject Classification:} 34C23, 34F10}

\section{Introduction} 		      

We study the one-dimensional nonlocal elliptic equation 
which comes from the logistic equation of population dynamics.

\begin{equation}
\left\{
\begin{array}{l}
-\left(\Vert u'\Vert_2^2+ \Vert u\Vert_{p+1}^{p+1} \right)^q u''(x)  + 
u(x)^p = \lambda u(x), \enskip 
x \in I:= (0,1),
\vspace{0.1cm}
\\
u(x) > 0, \enskip x\in I, 
\vspace{0.1cm}
\\
u(0) = u(1) = 0,
\end{array}
\right.
\end{equation}
where $p > 1$ and $q > 0$ are given constants satisfying 
$q < \frac{p-1}{p+1}$ and 
$\lambda > 0$ is a bifurcation parameter.

Nonlocal one-dimensional 
elliptic problems have been studied intensively. 
In particular, the main topics seem to be the study of the existence,  nonexistence and multiplicity of the solutions. We refer to [2-4, 7-14, 21] and the references therein. Besides, in the study of nonlinear elliptic problems, the bifurcation analysis 
is one of the popular concern and there are a lot of results which observe the nonlinear elliptic problems from a view point of bifurcation analysis. We refer to [1, 5, 6] and the references therein. 
Although the number of papers about bifurcation diagrams is smaller than that in nonlinear 
elliptic eigenvalue problems,
some research has been conducted recently in the 
field of nonlocal elliptic problems. 
We refer to [16-20, 22].

The objective here is to obtain the precise asymptotic formulas for the bifurcation curves of the equation (1.1) in $L^2$-framework. 
Precisely, $\lambda$ is parameterized by 
$\alpha = \Vert u_\lambda\Vert_2$ such as $\lambda(\alpha)$ and 
we study the asymptotic behavior of $\lambda(\alpha)$ as 
$\alpha \to \infty$

The equatin (1.1) is motivated by the standard nonlinear eigenvalue problem of 
logistic type. 
\begin{equation}
\left\{
\begin{array}{l}
-w''(x) + w(x)^p = \gamma w(x), \enskip 
x \in I,
\vspace{0.1cm}
\\w(x) > 0, \enskip x\in I, 
\vspace{0.1cm}
\\
w(0) = w(1) = 0.
\end{array}
\right.
\end{equation}
Let $\xi > 0$ be an arbitrary given constant. Then we know from [1, 15] that there exists a unique solution pair $(w_\xi, \gamma(\xi)) \in C^2(\bar{I}) \times \mathbb{R}_+$
of (1.2) with $\Vert w_\xi\Vert_2 = \xi$. Further, $\gamma$ is parameterized by $\xi$ such as $\gamma = \gamma(\xi)$, and 
it is called $L^2$-bifurcation curve. 

As far as the author knows, however, there are a few results 
concerning the precise global behavior of $\gamma(\xi)$. The reason 
seems to be that it is standard to treat the global shape of the bifurcation curve 
$\gamma$ of (1.2) in $L^\infty$-
framework from a view point of classical bifurcation theory. Indeed, in many cases, $\gamma$ is parameterized by $L^\infty$ norm 
of the solution $w_\gamma$ corresponding to $\gamma$, and 
it is represented as 
$\gamma = \gamma(\Vert w_\gamma \Vert_\infty)$.   
However, at the same time, it is equally important to treat the eigenvalue problem (1.2) 
in $L^2$-framework, since the global structure of 
$\gamma(\Vert w_\gamma \Vert_\infty)$ and $\gamma(\xi)$ are 
totally different from each other. We know (cf. [1]) that 
as $\Vert w_\gamma\Vert_\infty \to \infty$, 
\begin{eqnarray}
\gamma(\Vert w_\gamma \Vert_\infty) = \Vert w_\gamma\Vert_\infty^{p-1} + O(1).
\end{eqnarray} 
On the other hand, it was shown in [15, Theorem 1] that, for $\xi \gg 1$, 
\begin{eqnarray}
\gamma(\xi) &=& \xi^{p-1} + C_1\xi^{(p-1)/2} + \frac{1}{p-1}C_1^2 
+ o(1),
\end{eqnarray}
where
\begin{eqnarray}
C_1 = (p+3)\int_0^1 \sqrt{\frac{p-1}{p+1} - s^2 + \frac{2}{p+1}s^{p+1}}ds.
\end{eqnarray}
We briefly explain why the difference between (1.3) and (1.4) occurs. 
We see that (1.3) is affected only by the asymptotic 
behavior of $w_\gamma$ at the center of the interval $I$. On the other hand, since $w_\xi = \xi(1 + o(1))$ in the interior of $I$, $\gamma(\xi)$ is 
also affected by the asymptotic behavior of $w_\xi$ 
near the boundary of $I$. In fact, the second term of (1.4) comes from the slope of $w_\xi$ close to the boundary of $I$. 

Recently, motivated by (1.2), the nonlocal 
bifurcation problem related to (1.1) has been studied in [21]. 

\begin{equation}
\left\{
\begin{array}{l}
-\left(a_1\Vert v\Vert_q^2+ a_2\Vert v\Vert_2^2 \right) v''(x)  + 
v(x)^p = \mu v(x), \enskip 
x \in I = (0,1),
\vspace{0.1cm}
\\
v(x) > 0, \enskip x\in I, 
\vspace{0.1cm}
\\
v(0) = v(1) = 0.
\end{array}
\right.
\end{equation}
Motivated by (1.3) and (1.4), 
we studied (1.6) and following asymptotic formula has been obtained 
in [21].

\vspace{0.2cm}

\noindent
{\bf Theorem 1.0 [21].} {\it Assume that $p > 3$. Let 
$a_1, a_2 \ge 0$ be constants satisfying   
$a_1 + a_2 > 0$. Then for any given constant $\zeta > 0$, there exists a unique solution pair 
 $(v_\zeta, \mu(\zeta)) \in C^2(\bar{I}) \times \mathbb{R}_+$ of (1.6) satisfying $\Vert v_\zeta \Vert_2 = \zeta$. Furthermore, 
 as $\zeta \to \infty$, 
\begin{eqnarray}
\mu(\zeta) &=&\zeta^{p-1}\left\{1 + C_1(a_1+ a_2)^{1/2} \zeta^{-(p-3)/2}
+ O(\zeta^{-(p-3)})\right\}.
\end{eqnarray}
}

\vspace{0.2cm}

\noindent
We note that the difference betweewn (1.4) and (1.7) comes also from 
the difference of the slope of $w_\xi$ and $v_\zeta$. 

It should also 
be pointed out that the Kirchhoff function in (1.6) does not contain 
$\Vert v'\Vert_2$, which appears frequently in Kirchhoff 
type nonlocal problems. Motivated by this, 
we consider here the 
effect of the nonlocal term $\Vert u'\Vert_2$ to the asymptotics of the 
bifurcation curves and establish the asymptotic formulas for $\lambda(\alpha)$ as $\alpha \to \infty$, 
which are different from (1.4) and (1.7). By the present results, we understand well how the nonlocal term $\Vert u'\Vert_2$ 
in (1.1) gives effect to the 
asymptotic behavior of $\lambda(\alpha)$ as $\alpha \to \infty$. 

Now we state our results.

\vspace{0.2cm}

\noindent
{\bf Theorem 1.1.} {\it Assume that $p > 1$ and $0 < q < \frac{p-1}{p+1}$ 
are constants.   
Then there exists a constant $\alpha_0 > 0$ such that for $\alpha > \alpha_0$, there exists a 
unique solution pair 
 $(u_\alpha, \lambda(\alpha)) \in C^2(\bar{I}) \times \mathbb{R}_+$ of (1.1) satisfying $\Vert u_\alpha \Vert_2 = \alpha$, and as $\alpha \to \infty$, 
\begin{eqnarray}
\lambda(\alpha) &=&\alpha^{p-1}\left\{1 + C_1\alpha^{-\{p-1-q(p+1)\}/2} + o(\alpha^{\alpha^{-\{p-1-q(p+1)\}/2}})\right\}.
\end{eqnarray}
}

\vspace{0.2cm}

\noindent
{\bf Example 1.2.} Assume that $q = \frac{p-1}{2p}$. Then there exists a constant 
$\alpha_0 > 0$, which is obtained explicitly, such that for any $\alpha > \alpha_0$, the 
solution pair $(u_\alpha, \lambda(\alpha))$ 
is uniquely determined. Furthermore, as $\alpha \to \infty$, 
\begin{eqnarray}
\lambda(\alpha) &=&\alpha^{p-1}\left\{
1 + C_1\alpha^{-(p-1)^2/(4p)} + o(\alpha^{-(p-1)^2/(4p)})\right\}.
\end{eqnarray}

\vspace{0.2cm}

It should be mentioned that we do not know whether $\alpha_0$ obtained in Example 1.2 is optimal or not.

The crutial point of the 
proof of Theorem 1.1 is 
to show the the uniqueness of the solution pair  
$(u_\alpha, \lambda(\alpha))$ under 
the assumptions of the Kirchhoff function. 

The remainder of this paper is organized as follows. In Section 2, 
we prove Theorem 1.1 with the aid of the results in [15, 16], Taylor 
expansion and complicated direct calculation.

\section{Proof of Theorem 1.1} 		      

In what follows, we use the notations defined in Section 1. 
We begin with the existence of the solution pair 
$(u_\alpha, \lambda(\alpha))$.

\vspace{0.2cm}

\noindent
{\bf Lemma 2.1.} 
{\it Let $\alpha \gg 1$ be a fixed constant. Then there exists at least solution pair 
$(u_\alpha, \lambda(\alpha))$ of (1.1) satisfying $u_\alpha = h_\alpha 
w_{\xi_\alpha}$ for 
the constants 
$h_\alpha > 0$ and $\xi_\alpha > 0$.

}

\vspace{0.2cm}

\noindent
{\it Proof.} The road map of the proof is as follows. 

\noindent
{\it Step 1.} We find a solution of (1.1) of the form $u = h w_\xi$ for 
$\xi > \xi_0$, where $h > 0$ and $\xi_0 > 0$ are constants determined later. 

\noindent
{\it Step 2.} We show that $h$ is a function of $\xi$ and 
$h = h_\xi$ is a strictly increasing function 
for $\xi > \xi_0$. 

\noindent
{\it Step 3.} We put $\alpha = h_\xi\Vert w_\xi\Vert_2 = h_\xi\xi$. Then 
$\alpha$ is a strictly increasing function for $\xi > \xi_0$ and 
$\alpha \to \infty$ as $\xi \to \infty$. 
Then we find that for $\alpha > \alpha_0:= h_{\xi_0}\xi_0$, there 
exists a solution $u = u_\alpha$ for $\alpha > \alpha_0$. 
 
\noindent
{\it Proof of Step 1.} 
We look for $u_\alpha$, which is 
represented as $u_\alpha =h w_\xi$ for some 
$h > 0$ and $\xi > 0$. By (1.1), we have 
\begin{eqnarray}
-(h^2\Vert w_\xi'\Vert_2^2 + h^{p+1}
\Vert w_\xi\Vert_{p+1}^{p+1})^q 
hw_\xi'' + h^{p}w_\xi^p &=& \lambda(\alpha)h w_\xi.
\end{eqnarray}
We put $\beta:= 
(h^2\Vert w_\xi'\Vert_2^2 + h^{p+1}\Vert w_\xi\Vert_{p+1}^{p+1})^q$.  
By this and (2.1), we have 
\begin{eqnarray}
-\beta w_\xi'' + h^{p-1} w_\xi^p &=& \lambda(\alpha)w_\xi.
\end{eqnarray}
Namely, 
\begin{eqnarray}
-w_\xi''(x) + \frac{h^{p-1}}{\beta}w_\xi^p 
= \frac{\lambda(\alpha)}{\beta}w_\xi.
\end{eqnarray}
Now we consider the equation $h^{p-1} = \beta$, namely, 
\begin{eqnarray}
h^{p-1} = h^{2q}(\Vert w_\xi'\Vert_2^2 
+ h^{p-1}\Vert w_\xi\Vert_{p+1}^{p+1})^q.
\end{eqnarray}
This is equivalent to 
\begin{eqnarray}
h^{(p-1-2q)/q} = \Vert w_\xi'\Vert_2^2 
+ h^{p-1}\Vert w_\xi \Vert_{p+1}^{p+1}.
\end{eqnarray}
We put 
\begin{eqnarray}
D(\xi):= \Vert w_\xi'\Vert_2^2 + \frac{2}{p+1}\Vert w_\xi\Vert_{p+1}^{p+1}.
\end{eqnarray}
We know from [15, Lemma 2.2] that $D(\xi)$ is a strictly increasing function of 
$\xi > 0$. 
We put $y := h^{p-1} - \frac{2}{p+1}$. Then by (2.5), we have 
\begin{eqnarray}
\left(y + \frac{2}{p+1}
\right)^{(p-1-2q)/(q(p-1))} = \Vert w_\xi\Vert_{p+1}^{p+1}y + D(\xi).
\end{eqnarray}
Let $\xi_0 > 0$ satisfy $D(\xi_0) = \left(\frac{2}{p+1}
\right)^{(p-1-2q)/(q(p-1))}$. 
Since $q < \frac{p-1}{p+1}$, we see that $\frac{p-1-2q}{q(p-1)} > 1$. 
Then for $\xi > \xi_0$, we see that there exists a unique solution 
$y = y_\xi > 0$ 
satisfying (2.7). By this and putting 
$\lambda = \beta \gamma(\xi) 
= h_\xi^{p-1}\gamma(\xi)$ and $\alpha = h_\xi \xi$, we 
obtain the solution pair $(u_\alpha, \lambda(\alpha))$ 
of (1.1) with $\alpha = \alpha_\xi = h_\xi w_\xi$, where 
$h_\xi = (y_\xi + \frac{2}{p+1})^{1/(p-1)}$.

\noindent
{\it Proof of Step 2.} Let $\xi_0 < \xi_1 < \xi_2$. Since $D(\xi)$ is 
strictly increasing in $\xi > 0$, we 
have $D(\xi_0) < D(\xi_1) < D(\xi_2)$. Moreover, by [1, 15], we see that 
$w_{\xi_0(x)} < w_{\xi_1(x)} < w_{\xi_2(x)}$ for $x \in I$ and consequently, 
$\Vert w_{\xi_1}\Vert_{p+1}^{p+1} 
< \Vert w_{\xi_2}\Vert_{p+1}^{p+1}$. Therefore, for $\xi_0 < \xi_1 < \xi_2$, 
there exists a unique $y_1 < y_2$ satisfying (2.7). This implies that 
there exists a unique $h =h_{\xi_1}$ and $h_{\xi_2}$ satisfying 
(2.4) and $h_{\xi_1} < h_{\xi_2}$. 

\noindent
{\it Proof of Step 3.} 
By putting 
$\alpha_\xi = h_{\xi}\xi$, 
we see that $\alpha_\xi$ is strictly increasing for $\xi > \xi_0$. 
Further, $h_\xi \to \infty$ as $\xi \to \infty$. Indeed, if 
there exists a constant $\delta > 0$ and a sequence $\{\xi_\eta\} \subset 
\{\xi\}$, which is denoted by $\{\xi\}$ again,  
such that $h_{\xi} \le \delta$ for 
$\xi > \xi_0$, then we obtain a contradiction, since 
l.h.s of (2.5) is bounded, while r.h.s of (2.5) tends to $\infty$, since 
$\Vert w'_\xi\Vert_2^2 \sim \xi^{(p+3)/2}$ for $\xi \gg 1$ (cf. (2.16) 
below). 
By this, for a given $\alpha > \alpha_0$, there exists $\xi > \xi_0$ such 
that $\alpha = \alpha_\xi$ and $(u_\alpha, \lambda(\alpha)) 
= (h_\xi w_\xi, \beta\gamma(\xi))$ satisfies (1.1).   
Thus, the proof is complete. \qed

\vspace{0.2cm}

\noindent
{\bf Lemma 2.2.} {\it For $\alpha \gg 1$, the solution pair 
$(u_\alpha, \lambda(\alpha))$ of (1.1) is unique.} 

\noindent
{\it Proof.} Let $j = 1,2$. Assume that $(u_j, \lambda_j)$  satisfies (1.1) 
with $\Vert u_j \Vert_2 = \alpha$. We 
put $\beta_j:= (\Vert u_j'\Vert_2^2 
+ \Vert u_j\Vert_{p+1}^{p+1})^q$. By (1.1), we have 
\begin{eqnarray}
-\beta_j u_j'' + u_j^{p} = \lambda_j u_j.
\end{eqnarray}
We put $W_j:= h_j^{-1}u_j$, where $h_j := \beta_j^{1/(p-1)}$. 
By (2.8), we have 
\begin{eqnarray}
-W_j'' + \frac{h_j^{p-1}}{\beta_j}W_j^{p} = \frac{\lambda_j}
{\beta_j}W_j.
\end{eqnarray}
Then we have 
\begin{eqnarray}
-W_j'' + W_j^p = \frac{\lambda_j}
{\beta_j}W_j.
\end{eqnarray}
This implies that there exists a unique constant 
$\xi_j > 0$ for $j = 1, 2$ such that 
$(W_j, \frac{\lambda_j}
{\beta_j}) = (w_{\xi_j}, \gamma(\xi_j))$ satisfies (1.2).  
By repeating the argument in Lemma 2.1, we see that $h_jw_{\xi_j}$ 
satisfies (2.5) and (2.6). We show that 
if $\xi_1 < \xi_2$, then 
$h_1 < h_2$. To do this, as in the proof of Step 2 in Lemma 2.1 also 
holds in this case, and it is enough to show that $\xi_j \to \infty$ as $\alpha \to \infty$. 
Indeed, since $h_j\xi_j = \alpha \to \infty$, if $\xi_j \le C$, then $h_j \to \infty$. We know from (1.4) that $\gamma(\xi_j) \le C$, 
since $\xi_j \le C$, and 
by (1.2) and integration by parts, we have (2.5). 
This is a contradiction, since $(p-1-2q)/q > p-1$ and $h_j \to \infty$. 
Therefore, we see that $\xi_j \to \infty$ as $\alpha \to \infty$. 
Therefore, we see that the proof of Step 2 in Lemma 2.1 also holds and obtain that 
if $\xi_1 < \xi_2$, then 
$h_1 < h_2$. However, since $\alpha = h_1\xi_1 = h_2\xi_2$, we find 
$\xi_1 = \xi_2, h_1 = h_2$ and consequently, $w_{\xi_1} = w_{\xi_2}$. Namely, 
\begin{eqnarray}
u_1 = h_1w_{\xi_1} = h_2w_{\xi_2} = u_2, \enskip 
\lambda_1 = \beta_1\gamma(\xi_1) = h_1^{p-1}\gamma(\xi_1) = 
h_2^{p-1}\gamma(\xi_2) = \lambda_2.
\end{eqnarray} 
Thus, the proof of is complete. \qed

\vspace{0.2cm}

Now we give the proof of Example 1.2. In this case, we obtain $h_\xi$ 
explicitly. 

\vspace{0.2cm}

\noindent
{\bf Lemma 2.3.} {\it Let $q = (p-1)/(2q)$. Then there 
exists a constant $\alpha_0 > 0$, which is determined explicitly, such that 
for $\alpha > \alpha_0$, there exists a unique 
solution pair $(u_\alpha, \lambda(\alpha))$ of (1.1).   
}

\noindent
{\it Proof.} We see from the arguments in Lemmas 2.1 and 2.2, $\xi > \xi_0$ 
satisfying $\alpha = h_\xi \xi$ is unique if $h_\xi$ is strictly increasing 
for $\xi > \xi_0$. Since $q = \frac{p-1}{2p}$, from (2.5), we have 
\begin{eqnarray}
h^{2(p-1)} = \Vert w_\xi'\Vert_2^2 
+ h^{p-1}\Vert w_\xi \Vert_{p+1}^{p+1}.
\end{eqnarray}
By this, for a given $\xi > 0$, we have 
\begin{eqnarray}
h_\xi^{p-1} &=& \frac{\Vert w_\xi \Vert_{p+1}^{p+1} + 
\sqrt{\Vert w_\xi\Vert_{p+1}^{2(p+1)} + 4 \Vert w_\xi'\Vert_2^2}}
{2}
\\
&=& 
\frac{\Vert w_\xi \Vert_{p+1}^{p+1}+ 
\sqrt{4 (\Vert w_\xi'\Vert_2^2 
+ \frac{2}{p+1}\Vert w_\xi\Vert_{p+1}^{p+1}) 
+ \Vert w_\xi\Vert_{p+1}^{p+1}( \Vert w_\xi\Vert_{p+1}^{p+1} 
- \frac{8}{p+1})}}
{2}
\nonumber
\end{eqnarray}
We recall from (2.5) that $\Vert w_\xi' \Vert_{2}^{2} + 
\frac{2}{p+1}\Vert w_\xi\Vert_{p+1}^{p+1}$ is strictly increasing for $\xi > 0$. 
Let $\xi_0 > 0$ be a unique constant satisfying  
$\Vert w_{\xi_0}\Vert_{p+1}^{p+1} = \frac{8}{p+1}$. We know from [1] that 
$w_{\xi_1}(x) < w_{\xi_2}(x)$ ($0 < x < 1$) when $\xi_1 < \xi_2$. 
By this, we see that $\Vert w_\xi\Vert_{p+1}^{p+1}$ is srtictly 
increasing for $\xi > 0$. By this, (2.7) and (2.8), we see that $h_\xi$ 
is srtictly 
increasing for $\xi > \xi_0$. Further, $h_\xi \to \infty$ as $\xi \to \infty$. This implies that if $\alpha > \alpha_0 = h_{\xi_0}\xi_0$, then 
there exists a unique $\xi$ which safisfies $\alpha = h_\xi \xi$. 
Thus, the proof is complete. \qed

\vspace{0.2cm}

\noindent
{\bf Lemma 2.4.} {\it Assume that $h = h_\xi$ satisfies (2.5). Then as 
$\xi \to \infty$, 
\begin{eqnarray}
h_\xi = \xi^k + B\xi^{k - \frac{p-1}{2}}(1 + o(1)),
\end{eqnarray}
where $k = \frac{q(p+1)}{p-1-q(p+1)}$ and $B = \frac{k}{p+3}C_1$.
}

\noindent
{\it Proof.} We know from [15] that for $\xi \gg 1$, 
\begin{eqnarray}
\Vert w_\xi \Vert_{p+1}^{p+1} &=& \xi^{p+1} + \frac{p+1}{p+3}C_1 
\xi^{(p+3)/2} + o(\xi^{(p+3)/2}),
\\
\Vert w_\xi'\Vert_2^2 &=& \frac{2}{p+3}C_1\xi^{(p+3)/2} 
+ \frac{1}{p-1}C_1^2\xi^2 + o(\xi^2).
\end{eqnarray}
Indeed, by [15, (3.2) and (3.16)], we have (2.15). By this, (1.2) and (1.4), 
we have (2.16). Since $h_\xi \to \infty$ as $\xi \to \infty$, by (2.5), (2.15) and (2.16), 
we have 
$$
h_\xi^{(p-1-q(p+1))/q} = \xi^{p+1}(1 + o(1)).
$$
This implies that $h_\xi = \xi^k(1 + o(1))$. Now we put 
$h_\xi = \xi^k(1 + B(\xi))$. We write $B = B(\xi)$ for simplicity. 
We substitute $h_\xi$ into (2.5) and use (2.15) and (2.16) to obtain 
\begin{eqnarray}
&&\xi^{k(p-1-2q)/q}(1 + B)^{(p-1-2q)/q} 
\\
&&= \frac{2}{p+3}C_1\xi^{(p+3)/2} 
+ \frac{1}{p-1}C_1^2\xi^2 + o(\xi^2) 
\nonumber
\\
&&
\quad + \xi^{k(p-1)}(1 + B)^{p-1}(\xi^{p+1} + \frac{p+1}{p+3}C_1\xi^{(p+3)/2}
+ o(\xi^{(p+3)/2})).  
\nonumber
\end{eqnarray}
By this, Taylor expansion and direct calculation, we have 
\begin{eqnarray}
&&\xi^{p+1}(1 + B)^{(p-1-2q)/q} = \xi^{p+1}
\left(1 + \frac{p-1-2q}{q}B + o({\xi^{-m}})\right)
\\
&&= \xi^{-k(p-1)}\xi^{p+1}\left\{\frac{2}{p+3}C_1\xi^{-(p-1)/2} 
+ \frac{1}{p-1}C_1^2\xi^{-(p-1)} + o(\xi^{-(p-1)}\right\}
\nonumber
\\
&&+ 
\xi^{p+1}
\left\{
1 + (p-1)B+ \frac{p+1}{p+3}C_1\xi^{-(p-1)/2}(1 + o(1))
\right\}.
\nonumber
\end{eqnarray}
By this, we obtain 
\begin{eqnarray}
\frac{p-1-q(p+1)}{q}B = \frac{p+1}{p+3}C_1\xi^{-(p-1)/2}
(1 + o(1)).
\end{eqnarray}
By this, we see that $B = \frac{k}{p+3}C_1$. This implies (2.14). Thus, the proof is complete. \qed

\vspace{0.2cm}

\noindent
{\bf Lemma 2.5.} {\it As $\alpha \to \infty$, 
\begin{eqnarray}
\xi &=& \alpha^{1/(k+1)}\left\{1 - \frac{1}{k+1}B\alpha^{-(p-1)/(2(k+1))} 
+ o(\alpha^{-(p-1)/(2(k+1))})\right\}. 
\end{eqnarray}
}
\noindent
{\it Proof.} 
By Lemma 2.4, we have 
\begin{eqnarray}
\alpha &=& h_\xi \xi = \xi^{k+1}(1 + B\xi^{-(p-1)/2}(1 + o(1))).
\end{eqnarray}
By this, we have 
$\xi = \alpha^{1/(p+1)}(1 + o(1))$. By this and (2.21), we have 
\begin{eqnarray}
\xi^{k+1} &=& \alpha\{1 - B\xi^{-(p-1)/2}(1 + o(1))\}
\\
&=& \alpha\left\{1 - B\alpha^{-(p-1)/(2(k+1))}(1 + o(1))\right\}.
\nonumber
\end{eqnarray}
By this and Taylor expansion, we obtain (2.20). \qed

\vspace{0.2cm}

\noindent
{\it Proof of Theorem 1.1.} 
By (1.4), Lemma 2.4 and Taylor expansion, we have 
\begin{eqnarray}
\lambda &=& h^{p-1}\gamma(\xi) 
\\
&=& \xi^{k(p-1)}\{1 + (p-1)B\xi^{-(p-1)/2}(1 + o(1))\}
(\xi^{p-1} + C_1\xi^{(p-1)/2} + O(1))
\nonumber
\\
&=& \xi^{(k+1)(p-1)}\left\{1 + ((p-1)B + C_1)\xi^{-(p-1)/2} 
+ o(\xi^{-(p-1)/2})\right\}.
\nonumber
\end{eqnarray}
By this, Taylor expansion and Lemma 2.5, we have 
\begin{eqnarray}
\lambda &=& \alpha^{p-1}\left\{1 - \frac{1}{k+1}B\alpha^{-(p-1)/(2(k+1))}
(1 + o(1))
\right\}^{(k+1)(p-1)}
\\
&&\times \left\{1 + ((p-1)B + C_1)\alpha^{-(p-1)/(2(k+1))}(1 + o(1))\right\}
\nonumber
\\
&=& \alpha^{p-1}\left\{1 - (p-1)B\alpha^{-(p-1)/(2(k+1))}
(1 + o(1))
\right\}
\nonumber
\\
&&\times \left\{1 + ((p-1)B + C_1)\alpha^{-(p-1)/(2(k+1))}(1 + o(1))\right\}
\nonumber
\\
&=& \alpha^{p-1}\left\{1 + C_1\alpha^{-(p-1-q(p+1))/2}(1 + o(1))\right\}.
\nonumber
\end{eqnarray}
Thus, the proof is complete. \qed

\vspace{0.2cm}

\noindent
{\bf Example.} We consider the special case of 
Example 1.2. Let $p = 3$ and 
$q = (p-1)/(2p) = 1/3$. We know that $\lambda = h^2\gamma(\xi)$. 
Then by (2.13), (2.15), (2.16) and Taylor expansion, we have 
\begin{eqnarray}
h^2 &=& \frac{\Vert w_\xi\Vert_4^4 + \Vert w_\xi\Vert_4^4
\sqrt{1 + \frac{4\Vert w_\xi'\Vert_2^2}{\Vert w_\xi\Vert_4^8}}}{2}
\\
&=& \frac{\Vert w_\xi\Vert_4^4 + \Vert w_\xi\Vert_4^4\left\{1 + 
\frac{2\Vert w_\xi'\Vert_2^2}{\Vert w_\xi\Vert_4^8}(1 + o(1))\right\}}{2}
\nonumber
\\
&=& \Vert w_\xi\Vert_4^4 + \frac{\Vert w_\xi'\Vert_2^2}{\Vert w_\xi\Vert_4^4}(1 + o(1))
\nonumber
\\
&=& \xi^4 + \frac{2}{3}C_1\xi^3 + o(\xi^3).
\nonumber
\end{eqnarray}
By definition, we know $k = 2, B= \frac13 C_1$. Further, By Lemma 2.5, 
$\xi = \alpha^{1/3}(1 - \frac13 B\alpha^{-1/3} + o(\alpha^{-1/3}))$. 
By this, (1.4) and Taylor expansion, we obtain 
\begin{eqnarray}
\lambda &=& h^2\gamma(\xi) = (\xi^4 + \frac23C_1\xi^4 + o(\xi^4))
(\xi^2 + C_1\xi + O(1))
\\
&=& \xi^6 + \frac53C_1\xi^5 + o(\xi^5)
\nonumber
\\
&=& \alpha^2\{1 - \frac13B\alpha^{-1/3} + o(\alpha^{-1/3}\}^6 + \frac53C_1\alpha^{5/3}
+ o(\alpha^{5/3})
\nonumber
\\
&=& \alpha^2\{1 - 2B\alpha^{-1/3} + o(\alpha^{-1/3})\} + \frac53C_1\alpha^{5/3} + o(\alpha^{5/3})
\nonumber
\\
&=&\alpha^2\{1 - \frac23C_1\alpha^{-1/3} + o(\alpha^{-1/3})\}+ \frac53C_1\alpha^{5/3}
+ o(\alpha^{5/3})
\nonumber
\\
&=& \alpha^2 + C_1\alpha^{5/3} + o(\alpha^{5/3}).
\nonumber  
\end{eqnarray}

\end{document}